\documentclass[11pt]{amsart}

\title[The radius of analyticity for a problem in epitaxial growth]
{The radius of analyticity for solutions to a problem in epitaxial growth on the torus}
\author{David M. Ambrose}
\address{Department of Mathematics, Drexel University, 3141 Chestnut Street, Philadelphia, PA 19104, USA}
\email{dma68@drexel.edu}

\newtheorem{theorem}{Theorem}

\begin{document}

\begin{abstract}
A certain model for epitaxial film growth has recently attracted attention, with the existence of small global solutions
having being proved in both the case of the $n$-dimensional torus and free space.
We address a regularity question for these solutions, showing that in the case of the torus, the solutions
become analytic at any positive time, with the radius of analyticity growing linearly for all time.   
As other authors have, we take the Laplacian of the initial data to be in the Wiener algebra, and we find an explicit smallness
condition on the size of the data. Our particular condition on the torus is that the Laplacian of the initial data should have norm 
less than $1/4$ in the Wiener algebra.
\end{abstract}

\maketitle

\section{Introduction}

We study the equation
\begin{equation}\label{theEquation}
h_{t}=\Delta e^{-\Delta h},
\end{equation}
with spatial domain $\mathbb{T}^{n},$ the $n$-dimensional torus, 
subject to initial condition
\begin{equation}\label{initialCondition}
h(0,\cdot)=h_{0}.
\end{equation}
This equation has been derived in \cite{Krug1995}, and again more recently in \cite{marzuolaWeare}, as a model in 
epitaxial growth of thin films.  Two recent works have analyzed this model; Granero-Bellinchon and Magliocca have proved
global existence of small solutions on the torus, and decay to equilibrium \cite{graneroBellinchon}.  On free space rather
than the torus, Liu and Strain have demonstrated global existence of small solutions, decay to equilibrium, and 
analyticity of solutions \cite{liuStrain}.

Both the papers \cite{graneroBellinchon} and \cite{liuStrain} take initial conditions $h_{0}$ such that $\Delta h_{0}$ is small
in $\mathbb{A},$ the Wiener algebra.  The Wiener algebra on free space is the set of functions with Fourier transform in 
$L^{1},$ and on the torus is the set of functions with Fourier series in $\ell^{1}.$  In the present work, we also take 
initial data $h_{0}$ such that $\Delta h_{0}$ is small in the Wiener algebra, but we otherwise follow a different method
of working in spaces related to the Wiener algebra.  Like these works, we also find a smallness condition which is very explicit;
in particular, our existence theorem for small solutions on the torus will have the condition 
$\|\Delta h_{0}\|_{\mathbb{A}(\mathbb{T}^{n})}<1/4;$ using the same norm, Granero-Bellinchon and Magliocca require
the size of the data to be less than $1/10.$ 

We follow the approach of Duchon and Robert, who proved existence for all time of small vortex sheet solutions
with initial interface height small in the Wiener algebra \cite{duchonRobert}.  This method is to introduce a modification
of the Wiener algebra for functions on spacetime, with an exponential weight which implies analyticity of solutions at
positive times, then to make a fixed point formulation of the problem via Duhamel's formula, and then to get existence of
a fixed point from the contraction mapping theorem on these spaces.

Using these spaces as in \cite{duchonRobert}, one finds that the radius of analyticity of solutions grows linearly in time.
The author and Mazzucato have applied this technique previously to find analytic solutions of the two-dimensional
Kuramoto-Sivashinsky equations \cite{ambroseMazzucato}, and the author, Bona, and Milgrom have used it to find
analytic solutions of Bona-Chen-Saut systems \cite{ambroseBonaMilgrom}.  The author has also applied the technique
to find analytic solutions for mean field games \cite{ambroseMFG1}, \cite{ambroseMFG2}.

Other techniques to prove analyticity of solutions have come about from the fluid dynamics community, such
as the work of Grujic and Kukavica \cite{GK98}.  The method of Grujic and Kukavica was used by the author and
Mazzucato for the two-dimensional Kuramoto-Sivashinsky equation as well in \cite{ambroseMazzucato}, and the 
method results in a radius of analyticity which grows like $t^{1/2}$ for the Navier-Stokes equations and like $t^{1/4}$
for the Kuramoto-Sivashinsky equations.  The order of growth here is related to the order of the leading-order parabolic
operator in the evolution equations, which is second-order for Navier-Stokes and fourth-order for Kuramoto-Sivashinsky.
In the current problem on epitaxial growth \eqref{theEquation}, as will be seen below in \eqref{expanded}, the leading-order
parabolic term is again fourth-order.  The work of Liu and Strain \cite{liuStrain}
demonstrates analyticity of solutions, finding the radius of analyticity again growing like $t^{1/4}$ for sufficiently 
large times.

By proving the
analyticity of solutions of \eqref{theEquation} in the case that the spatial domain is the torus, and determining the
 linear-in-time increase in the radius of
analyticity, the present work complements the work of Granero-Bellinchon and Magliocca \cite{graneroBellinchon}.  
The present work also complements the work of Liu and Strain because while it is demonstrated in
\cite{liuStrain} that the radius of analyticity of solutions on free space grows like $t^{1/4}$ for sufficiently large times
(this is Theorem 4 of \cite{liuStrain}), and grows linearly at small times (see Proposition 11 of \cite{liuStrain}), we add a bit
of detail to this short-time result.  To be precise, we show that the linear growth rate of the radius of analyticity on free
space can be taken arbitrarily large, and that the time interval on which one has this growth can be taken arbitrarily large,
by making the smallness constraint on the data more stringent.  Moreover, we believe that the contrast between the Liu-Strain
result on the radius of analyticity on free space (growth like $t^{1/4}$) and the main theorem of the present work (Theorem 1
at the end of Section 2 below, which indicates linear growth for all time) is itself noteworthy.

We now discuss the formulation of the problem, which starts the same as the formulation in \cite{graneroBellinchon} and
\cite{liuStrain}, which is by making the Taylor expansion of the exponential in \eqref{theEquation}.
Making this Taylor expansion of the exponential, we can write \eqref{theEquation} as
\begin{equation}\nonumber
h_{t}=\sum_{j=0}^{\infty}\Delta\left(\frac{(-\Delta h)^{j}}{j!}\right).
\end{equation}
We rewrite this by separating out the first two terms from the sum, and introducing some notation:
\begin{equation}\label{expanded}
h_{t}=-\Delta^{2}h +\sum_{j=2}^{\infty}\Delta F_{j},
\end{equation}
where $F_{j}$ is given by
\begin{equation}\label{Aj}
F_{j}=\frac{(-1)^{j}}{j!}\left((\Delta h)^{j}\right).
\end{equation}
Integrating in time, we find a Duhamel formula for our solutions,
\begin{equation}\nonumber
h(t,\cdot)=e^{-\Delta^{2}t}h_{0}+\sum_{j=2}^{\infty}\int_{0}^{t}
e^{-\Delta^{2}(t-s)}\Delta F_{j}(s,\cdot)\ ds.
\end{equation}
We give the name $I^{+}$ to the integral operator appearing here, namely 
\begin{equation}\label{IPlus}
I^{+}f=\int_{0}^{t}e^{-\Delta^{2}(t-s)}\Delta f(s,\cdot)\ ds,
\end{equation}
and this gives the resulting form for the Duhamel integral
\begin{equation}\label{duhamel}
h(t,\cdot)=e^{-\Delta^{2}t}h_{0}+\sum_{j=2}^{\infty}I^{+}F_{j}.
\end{equation}
Following the method of Duchon and Robert, we view \eqref{duhamel} as a fixed point problem.
We note that since the mean of $h$ is conserved by \eqref{theEquation}, we may, without loss of generality,
assume that the mean of $h_{0}$ is equal to zero.

In Section 2, we prove the existence of global solutions in the case of the torus, introducing the Duchon-Robert-type
modifications of the Wiener algebra, demonstrating a bound for $I^{+},$ and giving our
contraction argument.  In Section 3, we make some remarks about how the method would work on 
$\mathbb{R}^{n}$ instead of the torus, with the most notable difference being that the result by this method is only then for a short time.

\section{Global solutions on the torus}

\subsection{Function spaces}\label{functionSpace}

We use function spaces based on the Wiener algebra, in which functions have absolutely summable Fourier series.
However, we both use polynomial and exponential weights, and take a version of the Wiener algebra for functions on
spacetime.  We define $\mathcal{B}_{\alpha}^{j},$ for $\alpha>0$ and $j\in\mathbb{N},$ to be the set of functions
on $[0,\infty)\times\mathbb{T}^{n},$ continuous in time, such that the following norm is finite:
\begin{equation}\label{definitionOfNorm}
\|f\|_{\mathcal{B}_{\alpha}^{j}}=\sum_{k\in\mathbb{Z}^{n}}|k|^{j}\sup_{t\in[0,\infty)}
e^{\alpha t |k|}|\hat{f}(t,k)|.
\end{equation}
We note that for any $\alpha$ and any $j,$ this is a Banach algebra.  We will only need this fact for
$\mathcal{B}_{\alpha}^{0},$ and we now establish this.  To begin, let $f\in\mathcal{B}_{\alpha}^{0}$ and
let $g\in\mathcal{B}_{\alpha}^{0}.$  To compute the norm of $fg,$ we begin plugging into \eqref{definitionOfNorm}:
\begin{equation}\nonumber
\|fg\|_{\mathcal{B}_{\alpha}^{0}}=\sum_{k\in\mathbb{Z}^{n}}\sup_{t\in[0,\infty)}e^{\alpha t |k|}
|\widehat{fg}(t,k)|.
\end{equation}
We next use the convolution formula for the transform of the product $fg,$ and we use the triangle inequality:
\begin{equation}\nonumber
\|fg\|_{\mathcal{B}_{\alpha}^{0}}\leq\sum_{k\in\mathbb{Z}^{n}}\sup_{t\in[0,\infty)}
e^{\alpha t |k|} \sum_{j\in\mathbb{Z}^{n}} |\hat{f}(t,k-j)||\hat{g}(t,j)|.
\end{equation}
We manipulate the supremum and the exponential:
\begin{multline}\nonumber
\|fg\|_{\mathcal{B}_{\alpha}^{0}}
\\
\leq \sum_{k\in\mathbb{Z}^{n}}\sum_{j\in\mathbb{Z}^{n}}
\left(\sup_{t\in[0,\infty)}\left(e^{\alpha t |k-j|}|\hat{f}(t,k-j)|\right)\right)
\left(\sup_{t\in[0,\infty)}\left(e^{\alpha t |j|}|\hat{g}(t,j)|
\right)\right).
\end{multline}
By summing first in $k$ and then in $j,$ we have our conclusion, namely
\begin{equation}\nonumber
\|fg\|_{\mathcal{B}_{\alpha}^{0}}\leq \|f\|_{\mathcal{B}_{\alpha}^{0}}\|g\|_{\mathcal{B}_{\alpha}^{0}}.
\end{equation}

\subsection{Operator estimate}\label{operatorBound}

We establish now a bound for the operator $I^{+}.$
Upon taking the Fourier transform of \eqref{IPlus}, we have
\begin{equation}\nonumber
\widehat{I^{+}f}(t,k)=-\int_{0}^{t}e^{-|k|^{4}(t-s)}|k|^{2}\hat{h}(s,k)\ ds.
\end{equation}
This operator, $I^{+},$ is bounded from the space $\mathcal{B}_{\alpha}^{0}$ to 
$\mathcal{B}_{\alpha}^{2},$ with this gain of derivatives because of the leading-order
parabolic term in the equation.  

As remarked upon above, we are only interested in functions with mean zero, so we 
let $f\in\mathcal{B}_{\alpha}^{0}$ have mean zero.  Then, we compute the norm of $I^{+}f$
in $\mathcal{B}_{\alpha}^{2}:$
\begin{equation}\nonumber
\|I^{+}f\|_{\mathcal{B}_{\alpha}^{2}}
=\sum_{k\in\mathbb{Z}^{n}\setminus\{0\}}|k|^{4}\sup_{t\in[0,\infty)}e^{\alpha t |k|}
\left|\int_{0}^{t}e^{-|k|^{4}(t-s)}\hat{f}(s,k)\ ds\right|.
\end{equation}
We use the triangle inequality and adjust factors of the exponential:
\begin{equation}\nonumber
\|I^{+}f\|_{\mathcal{B}_{\alpha}^{2}}\leq\sum_{k\in\mathbb{Z}^{n}\setminus\{0\}}|k|^{4}
\sup_{t\in[0,\infty)}e^{\alpha t |k|}\int_{0}^{t}e^{-|k|^{4}(t-s)}e^{-\alpha s|k|}e^{\alpha s|k|}|\hat{f}(s,k)|\ ds.
\end{equation}
We then make some manipulations with supremums:
\begin{multline}\nonumber
\|I^{+}f\|_{\mathcal{B}_{\alpha}^{2}}
\\
\leq
\sum_{k\in\mathbb{Z}^{n}\setminus\{0\}}|k|^{4}\sup_{t\in[0,\infty)}
e^{\alpha t |k|}
\int_{0}^{t}e^{-|k|^{4}(t-s)-\alpha s|k|}\left(
\sup_{\tau\in[0,\infty)}e^{\alpha \tau |k|}|\hat{f}(\tau,k)|\right)\ ds
\\
\leq\left(\sum_{k\in\mathbb{Z}^{n}\setminus\{0\}}\sup_{\tau\in[0,\infty)}e^{\alpha \tau |k|}|\hat{f}(\tau,k)|\right)
\\
\cdot
\left(\sup_{k\in\mathbb{Z}^{n}\setminus\{0\}}\sup_{t\in[0,\infty)}|k|^{4}e^{\alpha t|k|}
\int_{0}^{t}e^{-|k|^{4}(t-s)-\alpha s|k|}\ ds\right).
\end{multline}
The first factor on the right-hand side is just $\|f\|_{\mathcal{B}_{\alpha}^{0}},$ and we must compute the second
factor to ensure that it is finite.

We compute the integral:
\begin{equation}\nonumber
\int_{0}^{t}e^{|k|^{4}s-\alpha s|k|}\ ds
=\frac{e^{|k|^{4}t-\alpha t |k|}-1}{|k|^{4}-\alpha|k|}.
\end{equation}
The denominator here indicates a restriction on $\alpha;$ we take $\alpha\in(0,1),$ and the denominator is then
never equal to zero.
Using this, our relevant computation becomes
\begin{equation}\label{whereDiscretenessHelps}
\sup_{k\in\mathbb{Z}^{n}\setminus\{0\}}\sup_{t\in[0,\infty)}
\frac{(e^{|k|^{4}t-\alpha t |k|}-1)(e^{\alpha t |k|-|k|^{4} t})}{1-\alpha/|k|^{3}}.
\end{equation}
Our choice of $\alpha$ ensures that the denominator is positive, so the negative term in the numerator may be neglected.
We conclude
\begin{equation}\nonumber
\|I^{+}f\|_{\mathcal{B}_{\alpha}^{2}}\leq \left(\sup_{k\in\mathbb{Z}^{n}\setminus\{0\}}\frac{1}{1-\alpha/|k|^{3}}\right)
\|f\|_{\mathcal{B}_{\alpha}^{0}}
= \frac{1}{1-\alpha}\|f\|_{\mathcal{B}_{\alpha}^{0}}.
\end{equation}

\subsection{Contraction argument}\label{contraction}

We will prove existence of solutions in $\mathcal{B}_{\alpha}^{2};$ this means that we are requiring the Laplacian of 
our initial data, $\Delta h_{0},$ to be in the Wiener algebra.

We seek fixed points of the mapping $\mathcal{T},$ with $\mathcal{T}$ defined by the right-hand side of
\eqref{duhamel}:
\begin{equation}\nonumber
\mathcal{T}h=e^{-\Delta^{2}t}h_{0}+\sum_{j=2}^{\infty}I^{+}F_{j},
\end{equation}
where, of course, $F_{j}$ depends on $h$ through \eqref{Aj}.

We will show that $\mathcal{T}$ is a contraction on a closed ball in $\mathcal{B}_{\alpha}^{2}.$
For $h_{0}$ satisfying our stated condition, that $\Delta h_{0}$ is in the Wiener algebra,
it is straightforward to check that $e^{-\Delta^{2}t}h_{0}$ is in the space $\mathcal{B}_{\alpha}^{2}$
(and, in fact, we perform this calculation at the end of the section).
Furthermore, by taking $h_{0}$ small, we find that $e^{-\Delta^{2}t}h_{0}$ is small in $\mathcal{B}_{\alpha}^{2}.$
We let $X$ be the closed ball in $\mathcal{B}_{\alpha}^{2}$ centered at $e^{-\Delta^{2}t}h_{0},$ with radius $r_{1},$
with this radius to be determined.  We also let $r_{0}$ denote the norm of the center of our ball, i.e.,
\begin{equation}\nonumber
\|e^{-\Delta^{2} t}h_{0}\|_{\mathcal{B}_{\alpha}^{2}} = r_{0}.
\end{equation}
Then, note that we have
\begin{equation}\label{boundInBall}
\|f\|_{\mathcal{B}_{\alpha}^{2}}\leq r_{0}+r_{1},\qquad \forall f\in X.
\end{equation}

To show that $\mathcal{T}$ is a contraction on $X,$ we must show two things.  First, that $\mathcal{T}$ maps 
$X$ to $X.$  To this end, let $f\in X$ be given.  We must compute the norm of $\mathcal{T}f-e^{-\Delta^{2}t}h_{0}$
in $\mathcal{B}_{\alpha}^{2},$ and find that this norm is no more than $r_{1}.$  Using the triangle inequality as well as
 the operator bound of Section \ref{operatorBound}, we see that it is sufficient to show
\begin{equation}\nonumber
\frac{1}{1-\alpha}\sum_{j=2}^{\infty}\|F_{j}\|_{\mathcal{B}_{\alpha}^{0}}\leq r_{1}.
\end{equation}
Using the definition \eqref{Aj} and the algebra property of Section \ref{functionSpace}, we see that it is sufficient to show
\begin{equation}\nonumber
\frac{1}{1-\alpha}\sum_{j=2}^{\infty}\frac{(\|\Delta f\|_{\mathcal{B}_{\alpha}^{0}})^{j}}{j!}\leq r_{1}.
\end{equation}
We have $\|\Delta f\|_{\mathcal{B}_{\alpha}^{0}}=\|f\|_{\mathcal{B}_{\alpha}^{2}},$ and
also $\|f\|_{\mathcal{B}_{\alpha}^{2}}\leq r_{0}+r_{1}.$  So, if
\begin{equation}\label{mappingCondition}
\frac{1}{1-\alpha}\sum_{j=2}^{\infty}\frac{(r_{0}+r_{1})^{j}}{j!} \leq r_{1},
\end{equation}
then $\mathcal{T}$ maps $X$ to $X.$  This is one condition which $r_{0}$ and $r_{1}$ will need to satisfy.

By requiring the contracting property of $\mathcal{T},$ we will find another condition which must be satisfied.
We let $h$ and $\tilde{h}$ be in $X,$ and we compute the norm of the difference of $\mathcal{T}h$ and
$\mathcal{T}\tilde{h}:$
\begin{equation}\label{contractionDifferenceBound}
\|\mathcal{T}h-\mathcal{T}\tilde{h}\|_{\mathcal{B}_{\alpha}^{2}}\leq\frac{1}{1-\alpha}
\sum_{j=2}^{\infty}\|F_{j}-\tilde{F}_{j}\|_{\mathcal{B}_{\alpha}^{0}}.
\end{equation}
We may factor a difference of $j^{\mathrm{th}}$ powers as
\begin{equation}\nonumber
x_{1}^{j}-x_{2}^{j}=(x_{1}-x_{2})\sum_{\ell=0}^{j-1}x_{1}^{j-1-\ell}x_{2}^{\ell}.
\end{equation}
Using this with \eqref{Aj}, the algebra property for $\mathcal{B}_{\alpha}^{0},$ the triangle inequality, and
the bound \eqref{boundInBall}, we have the following estimate:
\begin{multline}\nonumber
\|F_{j}-\tilde{F}_{j}\|_{\mathcal{B}_{\alpha}^{0}}=
\left\|\frac{(-1)^{j}(\Delta h)^{j}}{j!}-\frac{(-1)^{j}(\Delta\tilde{h})^{j}}{j!}\right\|_{\mathcal{B}_{\alpha}^{0}}
=\frac{1}{j!}\left\|(\Delta h)^{j}-(\Delta\tilde{h})^{j}\right\|_{\mathcal{B}_{\alpha}^{0}}
\\
\leq\frac{1}{j!}\|\Delta h - \Delta\tilde{h}\|_{\mathcal{B}_{\alpha}^{0}}\left\|\sum_{\ell=0}^{j-1}
(\Delta h)^{j-1-\ell}(\Delta\tilde{h})^{\ell}\right\|_{\mathcal{B}_{\alpha}^{0}}
\\
\leq 
\frac{1}{j!}\|h-\tilde{h}\|_{\mathcal{B}_{\alpha}^{2}}\sum_{\ell=0}^{j-1}\|\Delta h\|_{\mathcal{B}_{\alpha}^{0}}^{j-1-\ell}
\|\Delta\tilde{h}\|_{\mathcal{B}_{\alpha}^{0}}^{\ell}
\leq \frac{1}{(j-1)!}\|h-\tilde{h}\|_{\mathcal{B}_{\alpha}^{2}}(r_{0}+r_{1})^{j-1}.
\end{multline}
Using this with \eqref{contractionDifferenceBound}, we find
\begin{equation}\nonumber
\|\mathcal{T}h-\mathcal{T}\tilde{h}\|_{\mathcal{B}_{\alpha}^{2}}\leq 
\left(\frac{1}{1-\alpha}\sum_{j=2}^{\infty}\frac{(r_{0}+r_{1})^{j-1}}{(j-1)!}\right)\|h-\tilde{h}\|_{\mathcal{B}_{\alpha}^{2}}.
\end{equation}
Recognizing a Taylor series, this becomes
\begin{equation}\nonumber
\|\mathcal{T}h-\mathcal{T}\tilde{h}\|_{\mathcal{B}_{\alpha}^{2}}\leq\frac{1}{1-\alpha}\left(e^{r_{0}+r_{1}}-1\right)
\|h-\tilde{h}\|_{\mathcal{B}_{\alpha}^{2}}.
\end{equation}
Our condition for $\mathcal{T}$ being a contraction, then, is that we need $r_{0}$ and $r_{1}$ to satisfy
\begin{equation}\label{contractingCondition}
\frac{e^{r_{0}+r_{1}}-1}{1-\alpha}<1.
\end{equation}

To ensure that the mapping $\mathcal{T}$ is a contraction, then, we must have $r_{0}$ and $r_{1}$ so that
\eqref{mappingCondition} and \eqref{contractingCondition} are satisfied.  
To ensure that these conditions are satisfied, we start by choosing $r_{1}=r_{0};$ then, we only need conditions on $r_{0}.$
With this choice, \eqref{contractingCondition} becomes 
\begin{equation}\label{exponentialCondition}
e^{2r_{0}}< 2-\alpha,
\end{equation}
or, upon taking the logarithm, 
\begin{equation}\nonumber
r_{0}<\frac{\ln(2-\alpha)}{2}.
\end{equation}

We next seek to satisfy \eqref{mappingCondition}.  Using the usual bound for the error in making a Taylor approximation, we
see that \eqref{mappingCondition} will be satisfied if
\begin{equation}\nonumber
\left(\frac{1}{1-\alpha}\right)\left(\frac{e^{2r_{0}}(2r_{0})^{2}}{2}\right)\leq r_{0}.
\end{equation}
In light of \eqref{exponentialCondition}, it is sufficient to choose $r_{0}$ such that
\begin{equation}\nonumber
\left(\frac{2-\alpha}{1-\alpha}\right)\left(2r_{0}\right)\leq 1.
\end{equation}
Thus, our mapping condition is
\begin{equation}\nonumber
r_{0}\leq\frac{1-\alpha}{2(2-\alpha)}.
\end{equation}.

It is possible to check using calculus that, for all $\alpha\in(0,1),$ 
\begin{equation}\nonumber
\frac{1-\alpha}{2(2-\alpha)}<\frac{\ln(2-\alpha)}{2}.
\end{equation}
Thus, our only constraint is $r_{0}\leq\frac{1-\alpha}{2(2-\alpha)}.$
As long as $r_{0}$ satisfies $r_{0}<1/4,$ there exists an $\alpha\in(0,1)$ such that
this constraint is satisfied.
As a final step, we wish to interpret our condition on $r_{0}$ for $\Delta h_{0}$ in the Wiener algebra rather than 
considering the spacetime function $e^{-\Delta^{2}t}h_{0}\in\mathcal{B}_{\alpha}^{2}.$
To this end, we compute the operator norm of the solution operator for the linear equation:
\begin{multline}\nonumber
\|e^{-\Delta^{2}t}h_{0}\|_{\mathcal{B}_{\alpha}^{2}}
=\sum_{k\in\mathbb{Z}^{n}}|k|^{2}\sup_{t\in[0,\infty)}e^{\alpha t|k|}e^{-|k|^{4}t}|\hat{h}_{0}(k)|
\\
\leq
\left(\sum_{k\in\mathbb{Z}^{n}}|k|^{2}|\hat{h}_{0}(k)|\right)
\left(\sup_{k\in\mathbb{Z}^{n}}\sup_{t\in[0,\infty)}e^{t|k|(\alpha-|k|^{3})}\right)
= \|\Delta h_{0}\|_{\mathbb{A}}.
\end{multline}
Thus, our condition on the data is $\|\Delta h_{0}\|_{\mathbb{A}} < 1/4.$  We have proved the following theorem.

\begin{theorem}
Let $h_{0}$ satisfy $\|\Delta h_{0}\|_{\mathbb{A}(\mathbb{T}^{n})}<1/4.$  Let $\alpha\in(0,1)$ be such that
\begin{equation}\nonumber
\|\Delta h_{0}\|_{\mathbb{A}}\leq \frac{1-\alpha}{2(2-\alpha)}.
\end{equation}
Then there exists $h\in\mathcal{B}_{\alpha}^{2}$ such that $h$ solves \eqref{theEquation} with initial data
\eqref{initialCondition}.  For any $t>0,$ the function $h(t,\cdot)$ is analytic, with radius of analyticity greater than
or equal to $\alpha t.$
\end{theorem}.

\section{Local solutions on $\mathbb{R}^{n}$}

While the free-space case was studied in detail in \cite{liuStrain}, we briefly mention now how 
the above method may be adapted to $\mathbb{R}^{n}$ instead of $\mathbb{T}^{n}.$ 
On the torus, we have taken advantage of the discreteness of the Fourier variable, specifically when  
estimating \eqref{whereDiscretenessHelps}.  With a continuous Fourier variable instead, the denominator in 
\eqref{whereDiscretenessHelps} can be negative, and in this case, we cannot estimate the second term in the
numerator for arbitrarily large $t.$  By restricting to a finite time interval, we can still prove existence of a solution
for which the radius of analyticity grows linearly in time.  We note that this restriction to a finite time interval
 is not always a requirement of the current method,
but instead depends somewhat on the equation under consideration.  As evidence for this, we note that the original
Duchon-Robert result for vortex sheets was on $\mathbb{R}$ rather than on a periodic interval, and 
the result was the existence of small global solutions \cite{duchonRobert}.

We do not repeat all of the details, since the only significant change is the estimate for the operator $I^{+}.$
Thus, we will give the definition of the function spaces, prove the operator estimate for $I^{+},$ and state a 
new theorem.

We denote our function spaces as 
$\mathcal{B}_{\alpha,T}^{j},$ where $\alpha>0$ again represent the rate of linear growth of the
radius of analyticity, $T>0$ is the length of the time interval considered, and $j\in\mathbb{N}$ is a Sobolev-type weight.
The functions in $\mathcal{B}_{\alpha,T}^{j}$ are continuous in time with the following norm being finite:
\begin{equation}\nonumber
\|f\|_{\mathcal{B}_{\alpha,T}^{j}}=\int_{\mathbb{R}^{n}}|\xi|^{j}\sup_{t\in[0,T]}e^{\alpha t|\xi|}|\hat{f}(t,\xi)|\ d\xi.
\end{equation}

We begin computing the $\mathcal{B}_{\alpha,T}^{2}$-norm of $I^{+}f,$ for some $f\in\mathcal{B}_{\alpha,T}^{0}:$
\begin{equation}
\|I^{+}f\|_{\mathcal{B}_{\alpha,T}^{2}}\leq \int_{\mathbb{R}^{n}}
|\xi|^{4}\sup_{t\in[0,T]}e^{\alpha t|\xi|}\int_{0}^{t}e^{-|\xi|^{4}(t-s)}|\hat{f}(s,\xi)|\ ds.
\end{equation}
We split the spatial integral into a piece over a compact set $\Omega\subseteq\mathbb{R}^{n}$ and its complement
$\mathbb{R}^{n}\setminus\Omega,$ where $\Omega$ is chosen such that for all $\xi\in\Omega,$ we have
\begin{equation}\nonumber
1-\frac{\alpha}{|\xi|^{3}}\leq \frac{1}{2}.
\end{equation}
That is to say, $\Omega$ is that set of all $\xi\in\mathbb{R}^{n}$ such that $|\xi|\leq (2\alpha)^{1/3}.$

We start by estimating the integral over $\Omega$ as follows:
\begin{multline}\nonumber
\int_{\Omega}|\xi|^{4}\sup_{t\in[0,T]}e^{\alpha t|\xi|}\int_{0}^{t}e^{-|\xi|^{4}(t-s)}|\hat{f}(s,\xi)|\ dsd\xi
\\
\leq C_{\alpha}\left(\int_{\Omega}\sup_{t\in[0,T]}e^{\alpha t|\xi|}|\hat{f}(t,\xi)|\ d\xi\right)
\left(\sup_{\xi\in\Omega}\sup_{t\in[0,T]}e^{\alpha t|\xi|}\int_{0}^{t}e^{-|\xi|^{4}(t-s)}e^{-\alpha s|\xi|}\ ds\right)
\\
\leq C_{\alpha,T}\|f\|_{\mathcal{B}_{\alpha,T}^{0}}.
\end{multline}
Here, we have used the compactness of the set $[0,T]\times\Omega$ to bound 
the factor of $|\xi|^{4}$ as well as
the double supremum on the right-hand side.

We next turn to the integral over $\mathbb{R}^{n}\setminus\Omega;$ our estimate for this case is very similar to our analysis
in Section \ref{operatorBound} above.  We find, similarly to the previous estimate, the following bound:
\begin{multline}\nonumber
\int_{\mathbb{R}^{n}\setminus\Omega}|\xi|^{4}\sup_{t\in[0,T]}e^{\alpha t |\xi|}\int_{0}^{t}
e^{-|\xi|^{4}(t-s)}|\hat{f}(s,\xi)|\ dsd\xi
\\
\leq \|f\|_{\mathcal{B}_{\alpha,T}^{0}}\left(\sup_{\xi\in\mathbb{R}^{n}\setminus\Omega}\sup_{t\in[0,T]}
|\xi|^{4}e^{\alpha t |\xi|}\int_{0}^{t}e^{-|\xi|^{4}(t-s)}e^{-\alpha s|\xi|}\ ds\right)
\\
\leq \left(\sup_{\xi\in\mathbb{R}^{n}\setminus\Omega}\frac{1}{1-\frac{\alpha}{|\xi|^{3}}}\right)
\|f\|_{\mathcal{B}_{\alpha,T}^{0}}\leq 2\|f\|_{\mathcal{B}_{\alpha,T}^{0}}.
\end{multline}
At this last step, we have used the definition of $\Omega$ to get the constant $2$ to appear.
We thus conclude that there exists $C_{\alpha,T}$ such that we have the operator norm bound
\begin{equation}\nonumber
\|I^{+}\|_{\mathcal{B}_{\alpha,T}^{0}\rightarrow\mathcal{B}_{\alpha,T}^{2}}\leq C_{\alpha,T}+2.
\end{equation}

Then, repeating the proof of Section \ref{contraction}, we have the following theorem.
\begin{theorem}  Let $\alpha>0$ and $T>0$ be given.  There exists $c>0$ such that for all $h_{0}$ satisfying 
$\|\Delta h\|_{\mathbb{A}(\mathbb{R}^{n})}<c,$ there exists a solution $h\in\mathcal{B}_{\alpha,T}^{2}$ to the problem
\eqref{theEquation}, \eqref{initialCondition}.  This solution is analytic at all times $t\in[0,T],$ with radius of analyticity greater
than or equal to $\alpha t.$
\end{theorem}

We remark that unlike in the case of the torus, we have not carefully tracked the size of the data for this theorem on 
free space.  We further remark that the linear rate of growth of the radius of analyticity, $\alpha,$ can be taken as large as
is desired, with the consequence that the amplitude threshold for existence on the time interval $[0,T]$ decreases to zero
as $\alpha$ increases to infinity.

\section*{Acknowledgments}
The author gratefully acknowledges support from the National Science Foundation through 
grant  DMS-1515849.  
\bibliography{ambroseEpitaxial}{}
\bibliographystyle{plain}

\end{document}